\documentclass{amsart}
\usepackage[OT4]{fontenc}
\usepackage{color}
\newtheorem{theorem}{Theorem}[section]
\newtheorem{lemma}[theorem]{Lemma}

\theoremstyle{definition}

\theoremstyle{remark}
\newtheorem{remark}[theorem]{Remark}

\numberwithin{equation}{section}

\def\Om{\Omega}
\def\C{\mathbb C^n}
\begin{document}

\title{On a problem of Chirka}

\author{S\l awomir Dinew}
\address{Faculty of Mathematics and Computer Science,  Jagiellonian University 30-348 Krakow, {\L}ojasiewicza 6, Poland}

\email{slawomir.dinew@im.uj.edu.pl}
\thanks{The first Author was supported by the Polish National Science Centre grant 2017/26/E/ST1/00955.}

\author{\.Zywomir Dinew}
\address{Faculty of Mathematics and Computer Science,  Jagiellonian University 30-348 Krakow, {\L}ojasiewicza 6, Poland}
\email{zywomir.dinew@im.uj.edu.pl}


\date{}


\begin{abstract}
We observe that a slight adjustment of a method of Caffarelli, Li, and Nirenberg yields that plurisubharmonic functions extend across subharmonic singularities as long as the singularities form a closed set of measure zero. This solves a problem posed by  Chirka. 
\end{abstract}

\maketitle

\section*{Introduction}
Let $\Om$ be a domain in $\C$ and $E\subseteq\Omega$ be a closed subset of Lebesgue measure zero. Let $u$ be a subharmonic function in $\Om$ which is furthermore plurisubharmonic in $\Om\setminus E$. Then $u$ is said to be a plurisubharmonic function with \textit{subharmonic singularities}
along $E$. It is natural to ask under what conditions on $E$ will $u$ be a genuine plurisubharmonic function on the whole $\Om$. Alternatively, one can ask when does the positive $(1,1)$-current $i\partial\overline{\partial} u$ on $\Om\setminus E$ extend past $E$ provided its 
trace measure $i\partial\overline{\partial} u\wedge \beta^{n-1}$, where $\beta:=i\partial\overline{\partial}\Vert z\Vert ^2$, extends past $E$ as a positive Borel measure. Such questions appear in potential theory on Teichm\"uller spaces- \cite{M} and 
in extension problems of positive line bundles
past small sets- \cite{S, Ha}. The problem is interesting only for sufficiently {\it large} sets $E$, as plurisubharmonic functions automatically extend past closed sets of 
$(2n-3)$-dimensional Hausdorff measure zero- \cite{Ha}.

In 2003,  Chirka  proved in \cite{Ch} the following result:  
\begin{theorem}
 Let $E$ be a real hypersurface of class $C^1$ in $\Om$. Then any function $u$ of the above type extends past $E$.
\end{theorem}
Under additional smoothness assumptions this was previously shown in \cite{B}, see also \cite{Ri} for a survey on these matters.

Note that the subharmonicity on $\Om$ is essential here as the example
$$u(z):=\begin{cases} \Vert z\Vert ^2& \text{ if }\ \Vert z\Vert \leq 1;\\
1\ {\rm}& \text{ if }\ \Vert z\Vert >1        \end{cases}
$$
shows. Unfortunately, as mentioned in \cite{Ch}, the method applied there crucially depends on the $C^1$ regularity. Thus in \cite{Ch} the following conjecture was made:

{\it Any plurisubharmonic function with subharmonic singularities along $E$ extends as a plurisubharmonic function provided $E$ is a closed subset of $\Om$ with a locally finite $(2n-1)$-dimensional Hausdorff measure.}

The aim of this note is to solve affirmatively this problem. In fact we observe that a slight modification of the tools introduced by Caffarelli, Li, and Nirenberg yields a much stronger result:

\begin{theorem}[Main Result]
Let $E\subseteq\Om$ be a closed subset of Lebesgue measure zero. Then any subharmonic function $u$ in $\Om$ which is plurisubharmonic in $\Om\setminus E$ is actually plurisubharmonic in the whole $\Om$. 
\end{theorem}

\section{Proof of the main result}
Of course there is nothing to prove in complex dimension $1$.

We begin with a classical result (see \cite{H}, Proposition 3.2.10') that subharmonicity and plurisubharmonicity can be tested through $C^2$-smooth local majorants. 
\begin{lemma}\label{visc}
Let $\Om\subseteq \C$ be a domain. An upper semicontinuous function $u$ on $\Om$ is subharmonic (respectively plurisubharmonic) if for every $z_0\in\Om$ and every local $C^2$-smooth function $\varphi$ defined near $z_0$
and satisfying $\varphi(z)\geq u(z)$ with equality at $z_0$ we have $\Delta\varphi(z_0)\geq 0$ (respectively we have $\frac{\partial ^2 \varphi}{\partial z_j\partial\bar{z}_k} (z_0)\geq 0$).
\end{lemma}
Note that $\varphi$ above is subject to conditions provided that it exists. In particular, the assumption can be void at some points $z_0\in\Om$. Also note that in the language of viscosity theory the above lemma says that $u$ is subharmonic iff it is a viscosity subsolution to the Laplace equation $\Delta v=0$ (that is $\Delta u\geq 0$) and $u$ is plurisubharmonic iff it is a viscosity subsolution to the {\it constrained complex Hessian} equation $\det^{+}\left(\frac{\partial ^2 v}{\partial z_j\partial\bar{z}_k}\right)= 0$ (that is 
$\det^{+}\left(\frac{\partial ^2 u}{\partial z_j\partial\bar{z}_k}\right)\geq 0$), where 
$${\det}^{+}(A)=\begin{cases}
              \det(A)\ & \text{ if }\ A\geq 0;\\
              -\infty\ & \text{ otherwise }.
             \end{cases}$$
For more details regarding the viscosity theory of such constrained complex Hessian equations we refer to \cite{EGZ}. 

The proof of the main result is essentially contained in Theorem 1.2 in \cite{CLN} once we adjust the tools applied there to the constrained complex equation 
$${\det}^{+}\left(\frac{\partial ^2 u}{\partial z_j\partial\bar{z}_k}\right)= 0.$$
We provide the details for the sake of completeness. We learned about this argument through \cite{BM}  and  \cite{Sw}.

As a direct corollary of Lemma \ref{visc} plurisubharmonicity of $u$ would follow if one can show that for any $z_0\in E$ and any local $C^2$ majorant $\varphi\geq u$, $\varphi(z_0)=u(z_0)$ one has 
$$\frac{\partial ^2 \varphi}{\partial z_j\partial\bar{z}_k} (z_0)\geq 0$$
as matrices.

Translating if necessary one may assume that $z_0$ is the origin, that $\Om$ contains a ball $B_{\delta_0}$ centered at the origin and that $\varphi$ is defined on $B_{\delta_0}$. For a fixed $0<\delta<\frac{\delta_0}2$ we consider the function

$$v_{\delta}(z):=\varphi(z)+\delta\Vert z\Vert ^2-\delta^3-u(z).$$
By the very definition of $\varphi$ we have $v_{\delta}(z)\geq 0$ on the collar $B_{2\delta}\setminus B_{\delta}$. Also, as $u$ is subharmonic, $v_{\delta}$ is bounded below on $B_{\delta}$ and $v_{\delta}$ is a lower semicontinuous viscosity supersolution to the Poisson equation
$\Delta v= \Delta\varphi+4n\delta=:f$, that is
$$\Delta v_{\delta}\leq f.$$
Note that $f$ is continuous.

Consider the following convex envelope in $B_{2\delta}$

$$\Gamma_{v_\delta}(z):=\sup\lbrace l(z)|\ l-\text{affine},\, l\leq v_{\delta}\, \text{ on }\  B_{\delta},\, l\leq 0\ \text{ on }\ B_{2\delta}\setminus B_{\delta} \rbrace.$$
As $v_\delta$ is bounded below the family of functions $l$ above is non void.

The following Alexandrov-Bakelman-Pucci type theorem for viscosity supersolutions (see Theorem 3.2 in \cite{CC}) is crucial:

\begin{lemma}
Let $v_\delta$ and $\Gamma_{v_{\delta}}$ be as above. Then for some universal constant $C$, dependent only on $n$ we have
$$\delta^3\leq \sup_{B_{\delta}}|v_{\delta}|\leq C\delta \left(\int_{\lbrace B_{\delta}\cap \lbrace v_{\delta}=\Gamma_{v_{\delta}}\rbrace\rbrace}\max\{f,0\}^{2n}dV\right)^{\frac1{2n}}.$$
\end{lemma}
We remark that this result is stated in \cite{CC} for continuous supersolutions but the  proof there works for lower semicontinuous supersolutions as well. Here we crucially exploit the fact that the Laplacian is a uniformly elliptic operator.

The upshot is that for every $\delta\in \left(0,\frac{\delta_0}2\right)$the function $v_\delta$ matches its convex envelope $\Gamma_{v_{\delta}}$ on a set of positive measure within $B_{\delta}$. As $E$ is of measure zero we pick a point 
$z_\delta\in\lbrace B_{\delta}\cap \lbrace v_{\delta}=\Gamma_{v_{\delta}}\rbrace\rbrace\setminus E$. As $u$ is plurisubharmonic around $z_{\delta}$ for any $0<r$ small enough and any unit vector $T\in\C$  we have
$$u(z_\delta)\leq\frac1{2\pi}\int_0^{2\pi}u(z_\delta+re^{i\theta}T)d\theta.$$
The same inequality is valid for the convex (hence plurisubharmonic) function $\Gamma_{v_\delta}$. But then 
$$\varphi(z_\delta)= u(z_\delta)+\delta^3-\delta\Vert z_\delta\Vert ^2+\Gamma_{v_{\delta}}(z_\delta)$$
$$\leq \delta^3-\delta\Vert z_\delta\Vert ^2+\frac1{2\pi}\int_0^{2\pi}[u+\Gamma_{v_\delta}](z_\delta+re^{i\theta}T)d\theta.$$
$$\leq \delta^3-\delta\Vert z_\delta\Vert ^2+\frac1{2\pi}\int_0^{2\pi}\left[\varphi(z_\delta+re^{i\theta}T)-\delta^3+\delta\Vert z_\delta+re^{i\theta}T\Vert ^2\right]d\theta$$
$$=\frac1{2\pi} \int_0^{2\pi}\varphi(z_\delta+re^{i\theta}T)d\theta+\delta r^2.$$
After dividing by $r^2$ and then letting $r\searrow 0^{+}$ we obtain
$$\sum_{j,k=1}^{n}\frac{\partial ^2\varphi}{\partial z_j\partial\bar{z}_k} (z_\delta)T_j\bar{T}_k\geq- \delta.$$
Finally, as $\delta\searrow 0^{+}$ we have $z_\delta\rightarrow 0(=z_0)$ and because $\varphi$ is $C^2$ it holds
$$\sum_{j,k=1}^{n}\frac{\partial ^2\varphi}{\partial z_j\partial\bar{z}_k} (z_0)T_j\bar{T}_k\geq 0.$$
As $T$ is an arbitrary unit vector this shows the non negative definiteness of the complex Hessian of $\varphi$ at $z_0$.
\begin{remark}Similar ideas are utilized in \cite{HL} to prove extension theorems but with the assumption of boundedness or continuity instead of subharmonicity. 
	\end{remark}
\begin{remark}
	It is well known that there exist closed sets of Lebesgue measure zero and of full Hausdorff dimension. Take for example a product of $2n$ copies of $A$, where $A=\{1\}\cup\bigcup_{j\geq 1} A_j$ and $A_{j}$ is a generalized Cantor set of Hausdorff dimension $1-1/j$ situated in the interval $[1-1/j,1-1/(j+1)]$. We see therefore that,  unlike in many similar removable singularity theorems, there is no Hausdorff dimension threshold up to which our theorem holds. Also often the possibility to extend objects is related to the vanishing of some capacity. Closed sets of measure zero can have positive capacity, as again the product of Cantor sets demonstrates. This is another distinctive feature of our theorem.
\end{remark}
\begin{remark} Our theorem could be combined with  some removable singularity theorems for (particular classes of) subharmonic functions, see e.g. \cite{Ga},\cite{Sa}. Now if  $u$ is subharmonic (of a particular class) on $\Omega\setminus F$, where $F$ is removable (for this particular class) and plurisubharmonic on $\Omega \setminus( E\cup F)$, where $E$ is closed and of measure zero, then $u$ is plurisubharmonic on $\Omega$.
\end{remark}
\begin{remark}\label{notclosed}
	We note that the closedness assumption on $E$ is not really necessary, and is introduced only  to ensure that  plurisubharmonicity on $\Omega\setminus E$ makes sense. Alternatively, we could assume that $E$ is any set of Lebesgue measure zero and if $u$ is plurisubharmonic on some neighborhood of $\Omega\setminus E$ then $u$ is plurisubharmonic on $\Omega$ if it is subharmonic there. This may be essential in some situations, as the closure of a null set can have positive measure.
	
\end{remark}

\bibliographystyle{amsplain}

\end{document}